\documentclass[11pt,a4paper]{article}
\usepackage[utf8]{inputenc}
\usepackage{amsfonts}
\usepackage{mathrsfs,amssymb,amsmath,wasysym} 
\usepackage{enumerate,textcomp}
\usepackage{amsthm}
\usepackage{pifont}
\usepackage{nameref}
\usepackage{booktabs}
\usepackage[numbers,square,sort&compress,comma]{natbib}
\usepackage{setspace}
\usepackage{color}
\usepackage{soul}
\usepackage[table]{xcolor}
\usepackage{graphicx}
\graphicspath{{../pdf/}{../jpeg/}{./image/}}
\DeclareGraphicsExtensions{.pdf,.jpeg,.png,.jpg,.eps}
\usepackage[left=2cm,top=1cm,right=2cm]{geometry}
\usepackage[T1]{fontenc}
\usepackage{hyperref}
\newtheorem*{theorem*}{Theorem}
\newtheorem{theorem}{Theorem}[section]
\newtheorem{lemma}[theorem]{Lemma}
\newtheorem{proposition}[theorem]{Proposition}

\newtheorem{assumption}{Assumption}
\numberwithin{equation}{section}

\allowdisplaybreaks
\title{\vspace{-0.in}\parbox{\linewidth }{\footnotesize\noindent}
Hybrid Ridgelet Deep Neural Networks for Data-Driven Arbitrage Strategies}

\author{Bahadur Yadav$^{}$\thanks{E-mail: bky0088@gmail.com}, Sanjay Kumar Mohanty $^{}$\thanks{Corresponding author. E-mail: sanjaymath@gmail.com} \\
	\footnotesize{Department of Mathematics, School of Advanced Sciences, Vellore Institute of Technology,}\\  \footnotesize{Vellore 632 014, Tamil Nadu, India}}
\date{ }

\begin{document}

\maketitle

\begin{abstract}
In this study, we propose a novel model framework that integrates deep neural networks with the Ridgelet Transform. The Ridgelet Transform on Borel measurable functions is used for arbitrage detection on high-dimensional sparse structures. This transform also enhances the expressive power of neural networks, enabling them to capture complex and high-dimensional market structures. Theoretically, we determine profitable trading strategies by optimizing hybrid ridgelet deep neural networks. Further, we emphasize the role of activation functions in ensuring stability and adaptability under uncertainty. We use a high-performance computing cluster for the detection of arbitrage across multiple assets, ensuring scalability, and processing large-scale financial data. Empirical results demonstrate strong profitability across diverse scenarios involving up to 50 assets, with particularly robust performance during periods of market volatility.
\end{abstract}


\section{Introduction}
In finance, arbitrage is the process of taking advantage of price variations for the same assets across marketplaces in an effort to make a profit without taking any risks. This is accomplished by simultaneously purchasing an asset at a lower cost and selling it at a higher price. Essentially, arbitrage is the practice of profiting from price differences by exploiting transient market inefficiencies. There are many different types of arbitrage, particularly in statistical and quantitative models.  Using short-term price differences for the same asset across different markets is known as classical, or pure, arbitrage (see \cite{cohen2020detecting}, \cite{cui2020detecting}). 
On the other hand, statistical arbitrage is a quantitative trading technique that finds and takes advantage of transient price differences between connected assets using mathematical models. It is based on the mean reversion concept, which states that prices eventually tend to revert to their historical averages. In order to make money, traders take advantage of these transient departures from anticipated relationships. \par
Pairs trading is a widely used form of statistical arbitrage that relies on the assumption that the spread between two strongly related assets is mean-reverting \cite{gatev2006pairs}. When the spread deviates from its historical equilibrium, traders expect it to revert over time, creating profit opportunities. However, the key drawback of this approach is its heavy reliance on the mean-reversion property. If this relationship breaks down and the spread fails to revert, the strategy may lead to losses rather than gains. We refer to \cite{rad2016profitability,acciaio2016model,burzoni2017model,burzoni2019pointwise,burzoni2021viability,davis2014arbitrage,fahim2016model,hobson2005static,neufeld2022deep,riedel2015financial}
 for theoretical concepts and developments in the literature on model-free arbitrage strategy. 
 \par
In contrast, Bondarenko\cite{bondarenko2003statistical} introduced a more general form of statistical arbitrage strategy, which is profitable on average, regardless of the final value of the underlying securities at maturity. Kassberger and Liebmann \cite{kassberger2017additive} expanded this idea by introducing G-arbitrage, which involves zero-cost payoffs. This generalization allows for more flexible trading strategies that can be adjusted based on available information. They linked the existence of G-arbitrage strategies to the absence of strategies fulfilling certain conditions, using G-measurable Radon-Nikodym densities. However, all of these approaches assume that the underlying securities follow a pre-determined probability measure.\par
Further, Lütkebohmert et al.\cite{lutkebohmert2021robust} introduced the concept of robust statistical arbitrage strategies, which refers to trading strategies that guarantee arbitrage profits under all probability measures within a given ambiguity set. This means the strategy remains profitable regardless of which measure from the set is the true one. It allows the strategy to succeed even in cases where traditional strategies like pairs trading would fail. However, their proposed method is based on linear programming, making it unsuitable for high-dimensional settings due to scalability limitations.\par
To address the limitations of linear programming, Neufeld et al.\cite{neufeld2024detecting} developed a deep neural network-based numerical method to identify robust statistical arbitrage strategies. Deep learning has been incorporated into robust statistical arbitrage more and more to address market inefficiencies in the face of uncertainty. Avellaneda and Lee established classical mean-reversion and cointegration models \cite{avellaneda2010statistical}, which were subsequently improved by machine learning techniques. While deep hedging models addressed frictions and uncertainty \cite{buehler2019deep}, deep neural networks demonstrated efficacy in equity prediction \cite{krauss2017deep}. Applications were further developed by reinforcement learning techniques, ranging from deep portfolio management \cite{deng2016deep} to recurrent models \cite{moody2001learning}. While studies bring attention to robustness challenges, distributionally robust hedging provides resilience to model and adversarial risks \cite{Zhang2020deep}, \cite{theate2023risk}. Combining deep learning and robust finance enables statistically profitable arbitrage strategies.
\par In trading, activation functions are crucial for modeling threshold-based behaviors, such as triggering buy signals or adjusting position sizes in response to signal strength\cite{apicella2021survey}.  In high-dimensional, reliable statistical arbitrage scenarios, functions like ReLU, sigmoid, and tanh efficiently capture these decision dynamics and guarantee smooth gradient flow, which is essential for stable deep network training.  Different activation functions have different characteristics. For example, some, like ReLU, are good at handling sparsity, while others, like tanh, capture saturation or probabilistic behavior.  Learning is improved with advanced functions such as GELU\cite{jagtap2023important}, Mish\cite{mathayo2022beta}, and SiLU\cite{elfwing2018sigmoid}, which offer smoother transitions and better gradient dynamics.
\par
Our work contributes to the literature on robust statistical arbitrage by focusing on the limitations of existing approaches, like the linear programming technique \cite{lutkebohmert2021robust}, and different neural network techniques, which are not scalable and profitable in diversified market scenarios. In order to address these issues, we present a recent Hybrid Ridgelet Deep Neural Network (HRDNN) technique for developing reliable statistical arbitrage strategies. The foundation of our strategy is a learning architecture that effectively approximates complex high-dimensional features essential to detecting arbitrage possibilities by fusing neural networks with the Ridgelet Transform. The Ridgelet Transform makes it easier to describe sparse information and provides a theoretical justification for neural networks' expressive capabilities.  Our methodology is resilient even in the absence of traditional statistical relationships like cointegration. Furthermore, we examine the ways in which activation functions affect the model's stability and capacity for generalization in the context of uncertainty, enhancing the dependability of the developed approaches in a range of market environments. Empirically, HRDNN consistently outperforms individual activations under a variety of transaction costs and is employed to detect both gradual market patterns and abrupt arbitrage opportunities. Using datasets with up to 50 underlying assets from the S\&P 500 universe and market stress conditions, empirical results show that our strategy is profitable and produces substantial profitability in a variety of scenarios. 
\par
The structure of the paper is as follows. Section \ref{2} introduces the underlying framework of the mathematical setting. In Section \ref{3}, we present our main theoretical findings, which are approximated using HRDNN. Section \ref{4} provides empirical evidence supporting the effectiveness of our method through a series of real-world examples. Finally, Section \ref{5} concludes the paper and outlines potential directions for future research.

\section{Mathematical Setting}\label{2}
Here, we analyze the price changes of a financial market with $a$ distinct assets at n\ future times, denoted by $t_1 < \dots < t_n$, where $a, n \in \mathbb{N}$. We assume that the price of each asset stays within a known and finite range. For each asset $j \in \{1, \dots, a\}$, we fixed two numbers  $\underline{U}^j$ and  $\overline{U}^j$ such that $\underline{U}^j < \overline{U}^j$, representing the lowest and highest values, where  $\underline{U}^j,  \overline{U}^j \in \mathbb{R}$. For each time period \( t_i \),  define  \( \Omega_i \) as
\[
\Omega_i = [\underline{U}^1, \overline{U}^1]^i \times \dots \times [\underline{U}^a, \overline{U}^a]^i \subset \mathbb{R}^{ia}, \quad i = 1, \dots, n ,
\]
this shows that the value of each asset is limited to its particular range at each time step up to time \( t_i \). 

Next, we define \( \Omega = \Omega_n \) that represents the region of every possible asset price path up to time \(t_n\). To model the evolution of asset prices over time, we introduce the process \( S = (S^{j}_{t_i})^{j=1,\dots,a}_ {i=1,\dots,n} \), which basically observes the price of each asset at every moment. At any time step $t_i$, the price of the asset $j$ is constrained by the bounds
\(\
S^j_{t_i} \in [\,\underline{U}^j, \overline{U}^j\,].
\)
 This process is defined as follows:
\[
S^{j}_{t_i} : \Omega \to [\underline{U}^j, \overline{U}^j]
\quad \text{such that}
\] 
\[
S^j_{t_i}\Big(x_1^1,x_2^1,\ldots,x_n^1;\;
               x_1^2,x_2^2,\ldots,x_n^2;\;
               \ldots;\
               x_1^a,x_2^a,\ldots,x_n^a\Big)
 = x_i^j, 
\]
where $x^j_i$ is the number of units for the asset $j$.

We now introduce a framework that characterizes the price of each asset at time \( t_i \). Let \( S_{t_i} = \left( S^j_{t_i} \right)^{j=1,\dots,a} \)
denote the price of the assets at time \( t_i \), where \(\
S_{t_i} = (S^1_{t_i}, S^2_{t_i}, \ldots, S^a_{t_i}) \in \mathbb{R}^a
\). At the initial time \(t_0\), the asset prices are deterministic, which means that they are known and fixed. Hence \( \left( S^j_{t_0} \right)^{j=1,\dots,a} \in [\underline{U}^1, \overline{U}^1] \times \cdots \times [\underline{U}^a, \overline{U}^a] \), and is called the observed spot price at time \(t_0\).

Next, we define a trading strategy  \( \Delta = \left( \Delta^j_i \right)^{j=1,\dots,a}_ {i=1,\dots,n-1} \), which contains a set of decision-making rules for every asset j at each time \(t_i\). Based on all the facts available at the time \(t_i\),  \( \Delta^j_i \), decide how much to invest in assets. For \( i \geq 1 \), each \( \Delta_i^j \) is a Borel measurable function and defined as \(\Delta_i^j : \mathbb{R}^{i  a} \to \mathbb{R} \)\quad \text{such that}
\[
\Delta_i^j(S_{t_1}, \ldots, S_{t_i}) = h_i^j ,
\]
where \(h_i^j\) is the trading position in asset \(j\) at time \(t_i\).
At the initial time \( t_0 \), the strategy \(\Delta_0^j \in \mathbb{R}\) is constant. By employing this strategy, we define the total gross profit from trading for each time step as \[
(\Delta \cdot S)_n = \sum_{j=1}^{a} \sum_{i=0}^{n-1} \Delta^j_i (S_{t_1}, \dots, S_{t_i}) (S^j_{t_{i+1}} - S^j_{t_i}),
\]
where \(\Delta^j_i (S_{t_1}, \dots, S_{t_i})\) represents the trading position of the assets at time \(t_i\), \(S^j_{t_{i+1}} - S^j_{t_i}\) is the change in the price of the assets and time \(t{_i}\).
Furthermore, we consider that the trading strategy \( \Delta \) leads to three different types of trading costs: transaction costs, borrowing costs associated with holding a short position, and liquidity costs due to bid-ask spreads. At any given time \(t_i\), the transaction cost occurs whenever a trader changes their asset position from one level \( \Delta^j_{i-1} \) to another \( \Delta^j_i \), this modification causes a transaction cost.
We discuss two particular forms of transaction cost that are frequently used in real-world scenarios \cite{buehler2019deep}, \cite{cheridito2017duality}.
\begin{enumerate}
    \item {Per-share transaction costs} (PSTC), which is defined by a non-negative measurable function 
\[
T^i_j : \mathbb{R} \times \mathbb{R} \to \mathbb{R}^+ \quad \text{such that} 
\]
    \[
    {T_i^ j}\bigl(S^j_{t_i}, x\bigr) = \lambda_T{_i^j} \ \ |x|,
    \]
    where \(\ \lambda_T{_i^j} > 0\) represents the fixed cost per unit traded  for asset $j$ at time \(t{_i}\).   
    
    \item {Proportional transaction costs} (PTC), which is defined by a non-negative measurable function 
\[
T^i_j : \mathbb{R} \times \mathbb{R} \to \mathbb{R}^+  \quad \text{such that} 
\]
     \[
    {T_i^ j}\bigl(S^j_{t_i}, x\bigr) =  \lambda_T{_i^j} \ \cdot S^j_{t_i}\ \cdot |x|,
    \]
    where \(\ \lambda_T{_i^j}  > 0\) represents the proportional cost rate for asset $j$ at time \(t{_i}\).
\end{enumerate}

In addition to transaction costs, liquidity costs arise from the bid-ask spread, which also penalizes large trades or agreements conducted when market liquidity is low. Considering bid-ask spreads are symmetric around the mid-price. These costs are defined by another non-negative measurable function 
\[
L_i^j : \mathbb{R} \times \mathbb{R} \to \mathbb{R}^+ \quad \text{such that} 
\]
\[
    {L_i^ j}\bigl(S^j_{t_i}, x\bigr)
=
0.5 \ \cdot \lambda_L{_i^j}  \ \cdot |x|,
\]

where \(\lambda_L{_i^j}  > 0\) represents the total amount of the bid-ask spread for asset $j$ at time \(t{_i}\).

Lastly, we include borrowing costs that are incurred when a trader enters a short position. This cost is define by a  non-negative measurable function \[ B_i^ j : \mathbb{R} \times \mathbb{R} \to \mathbb{R}^+ \quad \text{such that} \]
\[
    {B_i^ j}\bigl(S^j_{t_i}, x\bigr)
=
\lambda_B{_i^j} 
\ \cdot
\max\{-x, 0\}
\cdot
S^j_{t_i},
\]

where \(\lambda_B{_i^j}  > 0\) represents the {daily borrowing cost rate}, and a negative sign indicates that borrowing costs are only incurred when the trader has a short position.

Therefore, in order to determine the overall cost of carrying out a trading strategy, we add the costs incurred for each time step. The total cost is represented by \(E_n\), which is calculated as follows:
\[
E_n(\Delta) = \sum_{j=1}^{a} \sum_{i=0}^{n} 
\left[
T{_i^ j} \left( S^j_{t_i}, x \right)
+ L_i^ j \left( S^j_{t_i}, x \right) 
+ B{_i^ j} \left( S^j_{t_i}, x \right)
\right],
\]

and the total overall profit of a trading strategy \( \Delta \) is determined by \((\Delta \cdot S)_n -E_n(\Delta) \).

\section{Proposed Model}\label{3}
In this section, we discuss our key theoretical findings, which are framed in a generic way by an ambiguity set $\mathcal{P}$ of probability measures. The $\mathcal{P}$-Robust $\mathcal{G}$-Arbitrage strategy is presented, and the super replication problem that characterizes hedging costs under model uncertainty is also addressed.

\subsection{$\mathcal{P}$-Robust $\mathcal{G}$-Arbitrage Strategies}\label{3.1}
The idea of statistical arbitrage fundamentally depends on the characterization of an underlying probability measure. However, it is often difficult to determine which probabilistic model best captures the movements of asset prices.
To determine the model uncertainty, we consider a non-empty set \( \mathcal{P} \subseteq \mathcal{M}_1(\Omega) \) of Borel probability measures, representing admissible physical models. This ambiguity set captures uncertainty about the true underlying measure. We adopt the concept of a robust statistical arbitrage strategy and a multi-asset setting with transaction costs, as developed in \cite{lutkebohmert2021robust}.

Let \( \mathcal{P} \subseteq \mathcal{M}_1(\Omega) \), and \( \mathcal{G} \subseteq \sigma(S) \) be a sigma-algebra on \( \Omega \). A self-financing strategy \( \Delta \) is called a $\mathcal{P}$-Robust $\mathcal{G}$-Arbitrage strategy if the following conditions are satisfied:
\begin{enumerate}
    \item \( \mathbb{E}_P \left[ (\Delta \cdot S)_n -E_n(\Delta) \mid \mathcal{G} \right] \geq 0 \quad \text{P-a.s. for all} \ P \in \mathcal{P}, \)
    \item \( \mathbb{E}_P \left[ (\Delta \cdot S)_n -E_n(\Delta) \right] > 0 \quad \text{for some} \ P \in \mathcal{P}. \)
\end{enumerate}

Let \( B > 0 \), \( L > 0 \), and consider a set of strategies which are bounded by \( B \) and \( L \)-Lipschitz and defined as:
\begin{equation*}
W_{B,L} = \left\{ w_{c, \Delta} : \Omega \to \mathbb{R} \,\middle|\, 
w_{c, \Delta}(S) = c + (\Delta \cdot S)_n -E_n(\Delta) \right\},
\end{equation*}
for some \( c \in \mathbb{R} \),
\(
\Delta_i^j : \mathbb{R}^{ia} \to \mathbb{R} \text{ is } L\text{-Lipschitz} , \quad |c| \leq B, \\ 
\quad |\Delta_0^j| \leq B, \quad \text{and} \quad 
\|\Delta_i^j\|_{\infty, \Omega_i} = \sup_{x \in \Omega_i} |\Delta_i^j(x)| \le B.
\)

The constant \( B \) represents a budget constraint, which limits the maximum amount that an investor can allocate to their portfolio. The following lemma further establishes that the set \( {W}_{B,L}\) is compact. In order to ensure that the optimization problems are well-defined and the profit function is continuous with respect to trading strategies, we make the following assumption on the trading costs.
\begin{assumption}\label{assumption1}
We are assuming that all three types of trading costs are continuous.
\end{assumption}
\begin{lemma}\label{lemma1}(  \cite{neufeld2024detecting}):
Let Assumption \ref{assumption1} holds true. Then, for all \(B > 0\) and \(L > 0\), the space \( {W}_{B,L}\) is compact in the uniform topology on \(\Omega\).
\end{lemma}
Building on this, we consider the conditional super-replication problem. This allows us to characterize the minimal capital required to hedge a claim in a robust arbitrage framework.
 
Given a measurable function \( \Phi : \Omega \to \mathbb{R} \) and a sigma-algebra \( \mathcal{G} \subseteq \sigma(S) \), the conditional super-replication problem is defined as:
\[
\begin{aligned}
\mathcal{X}_{B,L}(\Phi,\mathcal{G}) = 
& \inf_{w_{c,\Delta}\in W_{B,L}} \{c \mid 
\mathbb{E}_P[w_{c,\Delta}(S)\mid\mathcal{G}]
 \geq &\mathbb{E}_P[\Phi(S)\mid\mathcal{G}],\;
\forall P\in\mathcal{P},\;\text{P-a.s.}\}
.\end{aligned}
\]

Solutions of  \(\mathcal{X}_{B,L}\) play an essential role in our framework because they capture robust arbitrage detection and robust pricing at the same time.

If the market does not accept any $\mathcal{P}$-robust $\mathcal{G}$-arbitrage strategy, $\mathcal{X}_{B,L}(\Phi, \mathcal{G})$ yields the highest price consistent with the absence of arbitrage for payoff $\Phi$, providing a robust upper bound on arbitrage-free pricing \cite{lutkebohmert2021robust}. In contrast, if the market does accept $\mathcal{P}$-robust $\mathcal{G}$-arbitrage strategy, then for the trivial claim $\Phi \equiv 0$, a value $\mathcal{X}_{B,L}(0, \mathcal{G}) < 0$ implies the existence of a strategy with negative initial cost and, conditional on $\mathcal{G}$, a non-negative payoff. In the following sections, we will approximate \(\mathcal{X}_{B,L}\) using the numerical approach as well as a ridgelet-based neural network. 

\subsection{An Approximation of \(\mathcal{X}_{B,L}\) through Penalization}\label{3.2}
In this section, we present a numerical approach to approximate  \( \mathcal{X}_{B,L} \) along with the corresponding \(\mathcal{P}\) -robust \(\mathcal{G}\)-arbitrage techniques. The strategy is based on penalization approaches, which are consistent with the strategies established in \cite{benamou2015iterative}, \cite{cominetti1994asymptotic}, \cite{eckstein2021robust}, \cite{eckstein2021computation}. We apply the assumption that the set \(\mathcal{P}\) of probability measures is finite in order to make the framework appropriate for numerical implementation and practical computation.
\begin{assumption}\label{assumption2}
The set \(\mathcal{P} \subset \mathcal{M}_1(\Omega)\) consists of finitely many probability measures.
\end{assumption}

Let \(M > 0\), and let \(\Phi : \Omega \to \mathbb{R}\) be a Borel-measurable payoff function. Let \(\mathcal{G} \subseteq \sigma(S)\) be a sub-\(\sigma\)-algebra of \(\sigma(S)\). To approximate \(\mathcal{X}_{B,L}(\Phi, \mathcal{G})\) numerically, we introduce the following penalized functional:
\begin{align*}
 \mathcal{X}_{B,L,k}(\Phi, \mathcal{G})=&\inf_{(\ w_{c, \Delta}) \in {W}_{B,L}} \Biggl\{ c + k \sum_{P \in \mathcal{P}} \int_\Omega \beta 
 \left( \mathbb{E}_P \left[\Phi(S) - {w_{c, \Delta} }(S) \mathcal{G} \right] \right) dP \Biggr\},  
\end{align*}
for each $k \in \mathbb{N}$, where \(\beta : \mathbb{R} \to \mathbb{R}^+\) is a penalization function designed to penalize trading strategies \(w_{c, \Delta} \in {W}_{B,L}\) that fail to satisfy the conditional super-replication constraint, 
\(\
\mathbb{E}_P[w_{c, \Delta}(S) \mid \mathcal{G}] \geq \mathbb{E}_P[\Phi(S) \mid \mathcal{G}] \quad \text{for all } P \in \mathcal{P}.
\)

To ensure that the penalization behaves correctly, we impose the following assumption on the geometry of \(\beta\).

\begin{assumption}\label{assumption3}
Let \(\beta : \mathbb{R} \to \mathbb{R}^+\) be a continuous function satisfying:
\[
\beta(x) = 0 \quad \text{if } x \leq 0, \qquad \beta(x) > 0 \quad \text{if } x > 0.
\]

\end{assumption}
As we will discuss in the following, Lemma \ref{lemma1} states that the compactness of \(W_{B,L}\) ensures the existence of a strategy \((w_{c*, \Delta*}) \in {W}_{B,L}\) that attains the minimum value of \(\mathcal{X}_{B,L,k}(\Phi, \mathcal{G})\) for each \(k\).

\begin{lemma}\label{lemma2}(  \cite{neufeld2024detecting})
 Let Assumptions \ref{assumption1}, \ref{assumption2}, \ref{assumption3} hold true. Let \( B > 0 \), \( L > 0 \), and let 
\(\Phi : \Omega \to \mathbb{R}\) be Borel measurable with 
\(\|\Phi\|_{\infty, \Omega} \leq B\). Then for every \(\sigma\)-algebra \( \mathcal{G} \subseteq \sigma(S) \) on \(\Omega\) and for all \( k \in \mathbb{N} \), the infimum of
\(\mathcal{X}_{B,L,k}(\Phi, \mathcal{G})\)
is attained, i.e., there exists \((w_{c*, \Delta*}) \in {W}_{B,L}\) 
{such that}
\begin{align*}
\mathcal{X}_{B,L,k}(\Phi, \mathcal{G}) =& c + k \sum_{P \in \mathcal{P}} \int_{\Omega} \beta\big(\mathbb{E}_P[\Phi(S)  - w_{c*, \Delta*}(S) \mid \mathcal{G}]\big) \, dP.
\end{align*}

\end{lemma}
Using Lemma \ref{lemma1} and Lemma \ref{lemma2}, we will show  the following result which demonstrates that
\(\mathcal{X}_{B,L}(\Phi, \mathcal{G})\) can  be approximated arbitrarily well by \(\mathcal{X}_{B,L,k}(\Phi, \mathcal{G})\)  by choosing k which is large enough.

\begin{proposition}\label{proposition1}(  \cite{neufeld2024detecting}) Let Assumptions \ref{assumption1}, \ref{assumption2}, \ref{assumption3} hold true. Let \( B > 0 \), \( L > 0 \), and let \( \Phi : \Omega \to \mathbb{R} \) be Borel measurable with \( \|\Phi\|_{\infty, \Omega} \leq B \). Then for every \( \sigma \)-algebra \( \mathcal{G} \subseteq \sigma(S) \), we have
\[
\lim_{k \to \infty} \mathcal{X}_{B,L,k}(\Phi, \mathcal{G}) = \mathcal{X}_{B,L}(\Phi, \mathcal{G}).
\]
\end{proposition}
In the following Section \ref{3.3}, we extend this approach by constructing an approximation of \(\mathcal{X}_{B,L}\) through a ridgelet-based neural network.

\subsection{An Approximation of \( \mathcal{X}_{B,L} \) using HRDNN}\label{3.3}
In this section, we introduce an alternative approach to approximate the functional \( \mathcal{X}_{B,L} \)  using HRDNN, which incorporate harmonic analysis techniques into trading formulation.

 We begin by introducing ridgelet functions as the fundamental building blocks of our approximation framework. For a fixed activation function \(\psi: \mathbb{R} \to \mathbb{R}\) and an input dimension m, a ridgelet function is defined as:
\[
R_{a,b,\psi}(x) = \psi(a \cdot x - b), \quad m \in \mathbb{N},
\]
where $ a \in \mathbb{R}^m$ represents the weight vector, $ b \in \mathbb{R}$ is the bias term,  and $ x \in \mathbb{R}^m$  is the input vector.
Intuitively, the ridgelet functions apply the activation function \(\psi\) to an affine translation of the input, just like neurons perform in standard neural networks.

Next, to build more expressive models, we introduce finite linear combinations of these ridgelet functions with arbitrary parameters. The class is defined as follows:
\[\mathrm{Ridge}_m = \bigcup_{N \in \mathbb{N}} \Bigl\{ f:\mathbb{R}^m \to \mathbb{R}^d \;\Bigm|\; f(x) = \sum_{j=1}^{N} \alpha_j \,\psi(a_j \cdot x - b_j) \Bigr\},\]
where
\(
\alpha_j \in \mathbb{R}^d,\quad 
a_j \in \mathbb{R}^m,\quad 
b_j \in \mathbb{R}
\).
Here, each function $f$ is composed of $N$ ridgelet terms, whose output weights mapping to the $ d$-dimensional target space are represented by \(\alpha_j\).

 To ensure that these functions behave in a well-controlled manner, we also apply Lipschitz continuity criteria and componentwise bounds. Specifically, for each domain  \(\Omega_i\), for \(i \in \mathbb{N}\), we define: 
\begin{equation*}
\begin{split}
\mathrm{Ridge}_{i,B,L} = \bigl\{ f \in \mathrm{Ridge}_m \;\big| \;
\pi_j \circ f|_{\Omega_i} \text{ is } L\text{-Lipschitz},
\|\pi_j \circ f\|_{\infty,\Omega_i} \leq B, \;\forall j=1,\ldots,a \bigr\},
\end{split}
\end{equation*}

where \(\pi_j : \mathbb{R}^a \to \mathbb{R}\) denotes the projection onto the \(j\)-th coordinate. In order to manage smoothness and growth within the given domains, it ensures that each output component is both bounded by a constant $B > 0$ and fulfils a Lipschitz condition with constant $L$.

By the universal approximation property of ridgelet systems,   Candès \cite{candes1999harmonic}, Eckstein \& Stephan\cite{eckstein2020lipschitz}, and Sonoda \cite{sonoda2017neural} ensure that every function in this class can be arbitrarily well approximated by HRDNN, making them an effective substitute for conventional neural network frameworks in the development of trading strategies. This property is crucial in our setting to approximate the trading strategies in \( {W}_{B,L} \), which require both boundedness and Lipschitz continuity. The following assumption is applied to the activation function to ensure that our neural networks preserve the universal approximation property. 
\begin{assumption}\label{assumption4}
We suppose that the activation function is either ReLU or a continuously differentiable, non-polynomial function.
\end{assumption}

\begin{lemma}[Density of Ridgelet Strategies]\label{lemma3}
Let Assumption \ref{assumption4} hold. Then, for each \(i = 1, \ldots, n\), \(j = 1, \ldots, a\), and for all \(B > 0\), \(L > 0\), the class \(\mathrm{Ridge}_{i,B,L}|_{\Omega_i}\) is dense in
\[
\left\{ f: \Omega_i \to \mathbb{R}^a \,\middle|\, \pi_j \circ f \text{ is } L\text{-Lipschitz}, \; \|\pi_j \circ f\|_{\infty,\Omega_i} \leq B \right\}
\]
with respect to the uniform topology on \(\Omega_i\).
\end{lemma}

{Proof of Lemma \ref{lemma3}}
Let \( B > 0 \), \( L > 0 \), \( 0 < \varepsilon < B \), and \( i \in \{1, \ldots, n\} \). Let
\(
f : \Omega_i \to \mathbb{R}^a
\)
be a function such that each component \( f_j = \pi_j \circ f \) is   $L$-Lipschitz and $\|f_j\|_{\infty, \Omega_i} \leq B$ for all $j = 1, \dots, a$. bounded by \( B \) on \(\Omega_i\). For each \( j = 1, \ldots, a \), define the truncated function  $\tilde{f}_j : \Omega_i \to \mathbb{R}$ by:
\[
\tilde{f}_j(x) = \max \left( \min \left( f_j(x), B - \frac{\varepsilon}{2 \sqrt{a}} \right), -B + \frac{\varepsilon}{2 \sqrt{a}} \right).
\]
This truncation preserves the Lipschitz constant, since for all \( x, y \in \Omega_i \), we have by construction
\[
|\tilde{f}_j(x) - \tilde{f}_j(y)| \leq |{f}_j(x) - {f}_j(y)| \leq L \|x - y\|.
\]
Then, according to the universal approximation theorem for Lipschitz functions  in
the form of [ Candès \cite{candes1999harmonic}, Eckstein \& Stephan\cite{eckstein2020lipschitz}, and Sonoda \cite{sonoda2017neural}, there exists for all \( j = 1, \ldots, a \) some ridgelet functions ${g}_j : \Omega_i \to \mathbb{R}$
\[
g_j(x) = \sum_{m=1}^N \alpha_j^m \psi_j^m(a_j^m \cdot x - b_j^m), \quad 
\]
{where each \( g_j \) is \( L \)-Lipschitz and bounded by \( B \) on \(\Omega_i\), and which fulfils}
\[
\|g_j - \tilde{f}_j\|_{\infty, \Omega_i} \leq \frac{\varepsilon}{2 \sqrt{a}}.
\]
Since by construction we have 
\(
\| \tilde{f}_j \|_{\infty, \Omega_i} \leq B - \frac{\varepsilon}{2\sqrt{a}}, \quad \text{it follows that for all \( j = 1, \ldots, a \),}
\)
\[\|g_j\|_{\infty,\Omega_i} \leq \|\tilde{f}_j\|_{\infty,\Omega_i} + \|g_j - \tilde{f}_j\|_{\infty,\Omega_i} \leq B - \frac{\varepsilon}{2 \sqrt{a}} + \frac{\varepsilon}{2 \sqrt{a}} = B.
\]

Furthermore, for the original function component \( f_j \),
\[
\|g_j - f_j\|_{\infty, \Omega_i} \leq \|g_j - \tilde{f}_j\|_{\infty, \Omega_i} + \|\tilde{f}_j - f_j\|_{\infty, \Omega_i} \leq \frac{\varepsilon}{2 \sqrt{a}} + \frac{\varepsilon}{2 \sqrt{a}} = \frac{\varepsilon}{\sqrt{d}}.
\]

Defining,  \(
g(x) = (g_1(x), \dots, g_a(x)) : \Omega_i \to \mathbb{R}^a,
\)
then \( g \in \text{Ridge}_m \), and since each \( g_j \) is \( L \)-Lipschitz and bounded by \( B \), we have \( g \in \text{Ridge}_{i,B,L} \). We get,
\[
\|g - f\|_{\infty, \Omega_i} \equiv \sup_{x \in \Omega_i} \sqrt{\sum_{j=1}^a |g_j(x) - f_j(x)|^2} \leq \sqrt{a \left( \frac{\varepsilon}{\sqrt{a}} \right)^2}\leq \varepsilon.
\]

Hence, the ridgelet-based function class \(\mathrm{Ridge}_{i,B,L}\) is dense in the set of all bounded, \(L\)-Lipschitz continuous functions on \(\Omega_i\).

This lemma ensures that any feasible trading strategy inside the class can be approximated by ridgelet-based networks. The restriction of the admissible strategy space \( {W}_{B,L}\) that can be represented by HRDNN leads to the set
\[
{W}^{\mathrm{Ridge}}_{B,L} = \left \{ w_{c, \Delta} \in {W}_{B,L} \,\middle|\, (\Delta^1_i, \Delta^2_i \ldots, \Delta^a_i) \in \mathrm{Ridge}_{i,B,L} \right\}.
\]
 \text{ for all } i = 1,\ldots,\(\ n-1\).
Accordingly, we define the ridgelet-based penalized functional:
\begin{align*}
   \mathcal{X}^{\mathrm{Ridge}}_{B,L,k}(\Phi,\mathcal{G}) =& \inf_{\ w_{c, \Delta} \in {W}^{\mathrm{Ridge}}_{B,L}} \Biggl\{ c + k \sum_{P \in \mathcal{P}} \int_\Omega \beta  \left( \mathbb{E}_P\left[\Phi(S) - {w_{c_, \Delta }}(S) \mid \mathcal{G} \right] \right) dP \Biggr\}. 
\end{align*}
Here, this method uses the same penalization framework  as \(\mathcal{X}_{B,L,k}(\Phi, \mathcal{G})\), which is implemented through the function \( \beta \). By Lemma \ref{lemma3}, we establish that \(\mathcal{X}^{\mathrm{Ridge}}_{B,L,k}(\Phi, \mathcal{G})\) approximates \(\mathcal{X}_{B,L,k}(\Phi, \mathcal{G})\) arbitrarily well by appropriate HRDNN in the following
 result.

\begin{proposition} \label{proposition2}
Let Assumptions \ref{assumption1}, \ref{assumption2}, \ref{assumption3}, \ref{assumption4} hold. Let \(L > 0\), \(B > 0\), and let \(\Phi : \Omega \to \mathbb{R}\) be a Borel-measurable payoff function satisfying \( \|\Phi\|_{\infty,\Omega} \leq B \). Then, for every sub-\(\sigma\)-algebra \(\mathcal{G} \subseteq \sigma(S)\) and for all \(k \in \mathbb{N}\), we have
\[
\mathcal{X}^{\mathrm{Ridge}}_{B,L,k}(\Phi, \mathcal{G}) = \mathcal{X}_{B,L,k}(\Phi, \mathcal{G}).
\]
\end{proposition}
{Proof of Proposition \ref{proposition2}.} Let \( \mathcal{G} \subseteq \sigma(S) \) be some $\sigma$
-algebra, and let \( k \in \mathbb{N} \). Since, we note that as \({W}^{\mathrm{Ridge}}_{B,L} \subseteq {W}_{B,L}\), we have
\(
\mathcal{X}^{\mathrm{Ridge}}_{B,L,k}(\Phi, \mathcal{G}) \geq \mathcal{X}_{B,L,k}(\Phi, \mathcal{G}).
\) Now we have to show that \(
\mathcal{X}^{\mathrm{Ridge}}_{B,L,k}(\Phi, \mathcal{G}) \leq \mathcal{X}_{B,L,k}(\Phi, \mathcal{G}).
\)

To prove the reverse inequality, we use  Lemma \ref{lemma2}, which ensures that there exists some \(w_{c, \Delta} \in {W}_{B,L}\) such that
\[
\mathcal{X}_{B,L,k}(\Phi, \mathcal{G}) = c + k \sum_{P \in \mathcal{P}} \int_{\Omega} \beta \left( \mathbb{E}_P[\Phi(S) - w_{c, \Delta}(S) \mid \mathcal{G}] \right) \, dP.
\] 

Next, by Lemma~\ref{lemma3}, for each \(i = 1, \ldots, n-1\), there exists a sequence of ridgelet functions 
\(
(\Delta^{1^{(m)}}_{i}, \ldots, \Delta^{d^{(m)}}_{i}) \in \mathrm{Ridge}_{i,B,L}
\)
such that
\[
(\Delta^{1^{(m)}}_{i}, \ldots, \Delta^{a^{(m)}}_{i}) \to (\Delta^1_i, \ldots, \Delta^a_i)
\quad \text{uniformly on } \Omega_i \text{ as } m \to \infty.
\]
Since both \(\Phi\) and the trading costs are bounded and continuous, and \(\Delta^{(m)}\) converges uniformly to \(\Delta\), we use the dominated convergence theorem along with the continuity of \(\beta\) to get
\small{
\begin{equation*}\nonumber
\begin{split}
    \lim_{m \to \infty} \left[ c + k \sum_{P \in \mathcal{P}} \int_{\Omega} \beta\left( \mathbb{E}_P[\Phi(S)-w_{c, \Delta^{(m)}}(S) | \mathcal{G}] \right) dP \right]
&= \left[ c + k \sum_{P \in \mathcal{P}} \int_{\Omega} \beta\left( \mathbb{E}_P[\Phi(S) - w_{c, \Delta}(S) |\mathcal{G}] \right)  dP \right] \\
&=\mathcal{X}_{B,L,k}(\Phi, \mathcal{G}).\\
\end{split}
\end{equation*}}

Thus, for each  \(w_{c, \Delta^{(m)}} \in {W}^{\mathrm{Ridge}}_{B,L}\), we obtain
\[
\mathcal{X}^{\mathrm{Ridge}}_{B,L,k}(\Phi, \mathcal{G}) 
\leq \lim_{m \to \infty} \left[ c + k \sum_{P \in \mathcal{P}} \int_{\Omega} \beta\left( \mathbb{E}_P[\Phi(S) - w_{c, \Delta^{(m)}}(S) \mid \mathcal{G}] \right) dP \right]
= \mathcal{X}_{B,L,k}(\Phi, \mathcal{G}).
\]

Combining both inequalities, we conclude that
\[
\mathcal{X}^{\mathrm{Ridge}}_{B,L,k}(\Phi, \mathcal{G}) = \mathcal{X}_{B,L,k}(\Phi, \mathcal{G}).
\]
This result ensures that HRDNN can be used to approximate the penalized super-replication functional \( \mathcal{X}_{B,L,k}(\Phi, \mathcal{G}) \). Building on this result, we now proceed from approximation theory to explicit computation. Following Section \ref{3.4} develops the computation of \(\mathcal{P}\)-robust statistical arbitrage strategies with HRDNN.

\subsection{Computation of \(\mathcal{P}\)-Robust Statistical Arbitrage Strategies using HRDNN} \label{3.4}
We focus on computing \(\mathcal{P}\)-robust statistical arbitrage strategies, which are a particular subclass of \(\mathcal{P}\)-robust \(\mathcal{G}\)-arbitrage strategies that can be obtained by choosing \( \mathcal{G} = \sigma(S_{t_n}) \). These strategies are theoretically and practically significant since they can be considered as profitable trading strategies that depend on market volatility. Our goal is to make an average profit, regardless of the specific terminal result. This is consistent with the idea of statistical arbitrage, which finds opportunities based on conditional expectations over the terminal value of the underlying asset (see \cite{bondarenko2003statistical}, \cite{lutkebohmert2021robust}, \cite{krauss2017statistical} ). However, it is extremely challenging to directly compute expectations conditioned on \( \sigma(S_{t_n}) \) from a numerical perspective. Because \( \sigma(S_{t_n}) \) is generated by a continuous random variable, it contains \emph{infinitely many measurable sets}, making direct computation of conditional expectations infeasible.

To address this challenge, we approximate \( \sigma(S_{t_n}) \) by employing a finite partition of the terminal state space. By discretizing the range of possible terminal values into an appropriate number of distinct domains.  Ridgelet Transformation allows us to formulate the problem in a finite-dimensional environment and makes the computation of conditional expectations feasible.

In order to solve this computational problem, Pseudo-Algorithm\cite{neufeld2024detecting} gives a systematic way to construct random sets from the sigma algebra \( \sigma(S_{t_n}) \). Specifically, for each \( i \in \mathbb{N} \), using the pseudo-algorithm, we construct the finite partition \( H_i \) of the terminal state space \([\underline{U}^1, \overline{U}^1] \times \dots \times [\underline{U}^a, \overline{U}^a]\). These partitions allow for feasible approximations by discretizing the state space into manageable subsets.
Furthermore, for each realization of Pseudo-Algorithm, we define:
\[
\mathcal{F}_i = \{ S_{t_n}^{-1}(A) : A \in \sigma(H_i) \}, \quad i \in \mathbb{N}.
\]

Here, each random set \( A_i \subset [\underline{U}^1, \overline{U}^1] \times \dots \times [\underline{U}^a, \overline{U}^a]\) is generated as a product of intervals of the following form:
\[
A_i = (p^{(i)}_1, q^{(i)}_1] \times \cdots \times (p^{(i)}_a, q^{(i)}_a]. \]

For each \( i \in \mathbb{N} \) and \( j = 1, \ldots, a \), let
\[
p^{(i)}_j \sim [\underline{U}^j, \overline{U}^j] , \quad q^{(i)}_j \equiv \overline{U}^j,
\]
where \( p^{(i)}_j \) a being independent of \( p^{(l)}_k \) for \((j,i) \neq (k,l)\). Moreover, we denote by \(\mathbb{P}^U\) the joint distribution of all 
\(
U = \{ (p^{(i)}_j, q^{(i)}_j) j=1,\ldots,a; \ i \in \mathbb{N} \} .
\)

The following proposition establishes an important result for almost all implementations of the Pseudo-Algorithm, the conditional expectations with respect to \( \sigma(S_{t_n}) \) can be arbitrarily well approximated by the conditional expectations determined with respect to \(\mathcal F_i \) provided that \( i \) is sufficiently large.
\begin{proposition} \label{proposition3}(  \cite{neufeld2024detecting}) Let \(L > 0\), \(B > 0\), and let \(\Phi : \Omega \to \mathbb{R}\) be Borel measurable with \(\|\Phi\|_{\infty,\Omega} \leq B\). Assume we have a trading strategy \(w_{{c, \Delta}} \in {W}_{B,L}\). Moreover, consider Pseudo-Algorithm, which generates random partitions using samples \(p^{(i)}_j, q^{(i)}_j\)  obtained from the distribution \(\mathbb{P}^U\)  for \(i \in \mathbb{N}\), \(j \in \{1, \dots, a\}\). Then we have:

(i) \(\sigma(S_{t_n}) = \sigma\left(\bigcup_{i=1}^\infty \mathcal{F}_i\right)\), \(\mathbb{P}^U\)-almost surely.

(ii) If \(\sigma(S_{t_n}) = \sigma\left(\bigcup_{i=1}^\infty \mathcal{F}_i\right)\), then for all \(\mathbb{P} \in \mathcal{P}\) we have
\[
\lim_{i \to \infty} \mathbb{E}_{\mathbb{P}} \left[ \Phi(S) - w_{c,\Delta}(S) \mid \mathcal{F}_i \right] 
= \mathbb{E}_{\mathbb{P}} \left[ \Phi(S) - w_{c,\Delta}(S) \mid \sigma(S_{t_n}) \right],
\]
  where the convergence holds both \(\mathbb{P}\)-almost surely and in \(L^1(\mathbb{P})\).

This result is essential since each \(\mathcal F_i \) is generated through a finite partition. As a result, it allows us to replace an unfeasible infinite-dimensional computation with an appropriate finite-dimensional approximation while maintaining the precision necessary for the formation of a robust strategy.

\end{proposition}
Building upon Proposition \ref{proposition1}, Proposition \ref{proposition2}, and Proposition \ref{proposition3}, we establish the following main theorem.

\begin{theorem*} \label{theorem}
Let \( B > 0 \), \( L > 0 \). We consider \(
\sigma(S_{t_n}) = \sigma\left( \bigcup_{i=1}^\infty \mathcal{F}_i \right)
\), this equality holds  true \( \mathbb{P}^{U} \)-a.s. as established by Proposition \ref{proposition3}. Then, by Lemma \ref{lemma2}, for all \( k \in \mathbb{N} \) and all \( i \in \mathbb{N} \), there exists a subsequence of the strategies \( w_{c^{(i)}, \Delta^{(i)}} \in W_{B,L} \) such that
\[
\mathcal{X}_{B,L,k}(\Phi, \mathcal{F}_i) = c^{(i)} + k \sup_{\mathbb{P} \in \mathcal{P}} \int_{\Omega} \beta \left( \mathbb{E}_P[\Phi(S) - w_{c^{(i)}, \Delta^{(i)}}(S) \mid \mathcal{F}_i] \right) \, d\mathbb{P},
\]
and some \( w_{c_{\infty}, \Delta_{\infty}} \in {W}_{B,L} \) such that
\[
\mathcal{X}_{B,L,k}(\Phi, \sigma(S_{t_n})) = c_{\infty} + k \sup_{\mathbb{P} \in \mathcal{P}} \int_{\Omega} \beta \, (\mathbb{E}_\mathbb{P} \left[ \Phi(S) - w_{c_{\infty}, \Delta_{\infty}}(S) \mid \sigma(S_{t_n})] \right) \, d\mathbb{P}.
\]
According to Lemma \ref{lemma1}, the set \( W_{B,L} \) is compact. Thus, there exists a subsequence of \( (w_{c^{(i)}, \Delta^{(i)}})_{i \in \mathbb{N}} \) such that
\(
(w_{c^{(i)}, \Delta^{(i)}})_{i \in \mathbb{N}} \to w_{c, \Delta} \in W_{B,L} \quad \text{uniformly as \( i \to \infty \).}
\)
Then, using Proposition \ref{proposition3}, the dominated convergence theorem, and the continuity of the trading cost functions and the penalty function \( \beta \), we can determine that

\begin{equation}\label{2.3}
  \begin{split}
    \mathcal{X}_{B,L,k}(\Phi, \sigma(S_{t_n})) 
&\leq c + k \sup_{\mathbb{P} \in \mathcal{P}} \int_{\Omega} \beta \, (\mathbb{E}_\mathbb{P} \left[ \Phi(S) - w_{c,\Delta}(S) \mid \sigma(S_{t_n}) \right]) \, d\mathbb{P}\\
&= \lim_{i \to \infty} \left( c + k \sup_{\mathbb{P} \in \mathcal{P}} \int_{\Omega} \beta \, (\mathbb{E}_\mathbb{P} \left[ \Phi(S) - w_{c,\Delta}(S) \mid \mathcal{F}_i \right]) \, d\mathbb{P} \right)\\
&= \lim_{i \to \infty} \lim_{j \to \infty} \left( c^{(j)} + k \sup_{\mathbb{P} \in \mathcal{P}} \int_{\Omega} \beta \, (\mathbb{E}_\mathbb{P} \left[ \Phi(S) - w_{c^{(j)},\Delta^{(j)}}(S) \mid \mathcal{F}_i \right]) \, d\mathbb{P} \right)  
  \end{split}  
\end{equation}
By applying the tower property of conditional expectation, along with Jensen's inequality, since \( \beta \) is convex, it holds for every \( j \geq i \), using \( \mathcal{F}_i \subseteq \mathcal{F}_j \), that
\begin{equation*}\label{2.4}
    \begin{split}
        \int_{\Omega} \beta \left( \mathbb{E}_\mathbb{P} \left[ \Phi(S) - w_{c^{(j)}, \Delta^{(j)}}(S) \mid \mathcal{F}_i \right] \right) \, d\mathbb{P}&=\int_{\Omega} \beta \left( \mathbb{E}_\mathbb{P} \left[ \mathbb{E}_\mathbb{P} \left[ \Phi(S) - w_{c^{(j)}, \Delta^{(j)}}(S) \mid \mathcal{F}_j \right] \mid \mathcal{F}_i \right] \right) \, d\mathbb{P}\\
&\leq  \int_{\Omega} \mathbb{E}_\mathbb{P} \left[ \beta \left( \mathbb{E}_\mathbb{P} \left[ \Phi(S) - w_{c^{(j)}, \Delta^{(j)}}(S) \mid \mathcal{F}_j \right] \right) \mid \mathcal{F}_i \right] \, d\mathbb{P}\\
&=\int_{\Omega} \beta \left( \mathbb{E}_\mathbb{P} \left[ \Phi(S) - w_{c^{(j)}, \Delta^{(j)}}(S) \mid \mathcal{F}_j \right] \right) \, d\mathbb{P}
    \end{split}
\end{equation*}
Thus, we obtain by (\ref{2.3}) that
\begin{equation*}
    \begin{split}
        \mathcal{X}_{B,L,k}(\Phi, \sigma(S_{t_n})) 
&\leq \lim_{j \to \infty} \left( c^{(j)} + k \sup_{\mathbb{P} \in \mathcal{P}} \int_{\Omega} \beta \left( \mathbb{E}_\mathbb{P} \left[ \Phi(S) - w_{c^{(j)}, \Delta^{(j)}}(S) \mid \mathcal{F}_j \right] \right) \, d\mathbb{P} \right)\\
&= \lim_{j \to \infty} \mathcal{X}_{B,L,k}(\Phi, \mathcal{F}_j)\\
&  \leq \lim_{j \to \infty} \left( c_{\infty} + k \sup_{\mathbb{P} \in \mathcal{P}} \int_{\Omega} \beta \left( \mathbb{E}_\mathbb{P} \left[ \Phi(S) - w_{c^{(j)}, \Delta^{(j)}}(S) \mid \mathcal{F}_j \right] \right) \, d\mathbb{P} \right)\\
&= c_{\infty} + k \sup_{\mathbb{P} \in \mathcal{P}} \int_{\Omega} \beta \, (\mathbb{E}_\mathbb{P} \left[ \Phi(S) - w_{c_{\infty}, \Delta_{\infty}}(S) \mid \sigma(S_{t_n}) \right]) \, d\mathbb{P}\\
&= \mathcal{X}_{B,L,k}(\Phi, \sigma(S_{t_n})).
    \end{split}
\end{equation*}

This means we have shown that for each \( k \in \mathbb{N} \),
\[
\lim_{j \to \infty} \mathcal{X}_{B,L,k}(\Phi, \mathcal{F}_j) = \mathcal{X}_{B,L,k}(\Phi, \sigma(S_{t_n})).
\]

Therefore, we conclude with Proposition \ref{proposition1} and Proposition \ref{proposition2} that
\[
\lim_{k \to \infty} \lim_{i \to \infty} \mathcal{X}^{\mathrm{Ridge}}_{B,L,k}(\Phi, \mathcal{G})
=
\lim_{k \to \infty} \lim_{i \to \infty} \mathcal{X}_{B,L,k}(\Phi, \mathcal{F}_i)
=
\lim_{k \to \infty} \mathcal{X}_{B,L,k}(\Phi, \sigma(S_{t_n}))
=
\mathcal{X}_{B,L}(\Phi, \sigma(S_{t_n})).
\]
\end{theorem*}
In particular, this theorem shows that for large enough \( k \) and \( i \), the robust value $\mathcal{X}_{B,L}(\Phi, \sigma(S_{t_n}))$ can be well approximated by $\mathcal{X}^{\mathrm{Ridge}}_{B,L,k}(\Phi, \mathcal{F}_i)$. Here, $\mathcal{X}^{\mathrm{Ridge}}_{B,L,k}(\Phi, \mathcal{F}_i)$ represents the optimization problem that is resolved by HRDNN, which provides the desired approximation precision while being computationally efficient.

In the following section, we will show the practical utility of the HRDNN approximation approach using real diversified financial data to show that our model is better compared to the existing model.

\section{Empirical Evaluation on Real-World Financial Data}\label{4}
This section presents the practical implementation of our proposed method on real-world financial data. The results provide empirical evidence supporting the applicability of our framework to physical measurements, even in a multi-asset context. Notably, our robust trading strategy demonstrates strong performance, particularly during periods of financial crisis.

The following sections detail the design and evaluation of our neural network-based arbitrage approach. Section~\ref{4.1} introduces the network architecture employed within our framework. Sections~\ref{4.2} and~\ref{4.3} examine the performance metrics and ratios across various activation functions. A numerical comparison with a benchmark model is discussed in Section~\ref{4.4}. Lastly, Section~\ref{4.5} investigates how different activation functions influence strategy performance in high-dimensional asset settings, and we will see that our model is superior to the existing model.

\subsection{Neural Network Architecture and Training Setup}\label{4.1}
We train our neural networks using the PyTorch framework, incorporating a batch normalization layer at the input stage. The architecture comprises three fully connected layers with \(32a\), \(64a\), and \(128a\) neurons, respectively, where \(a\) represents the number of assets. We employed a learning rate of $1e^{-4}$ for high-dimensional circumstances (10 assets or more). Empirical experiments have demonstrated that applying appropriate learning rates for high-dimensional scenarios, respectively, significantly enhances the training speed. We consider zero transaction costs (ZTC), PTC with \(\ \lambda_T{_i^j}  =0.001\), PSTC  with \(\ \lambda_T{_i^j}  =0.01\)  bid-ask spreads with \(\lambda_L{_i^j}  = 0.0002\), and daily borrowing costs with \(\lambda_B{_i^j}  =  0.1/252\). We test a variety of activation function types across all hidden layers to evaluate their impact on model performance.  Our results show that the model's performance changes noticeably depending on the chosen activation function. These effects are illustrated and analyzed in Tables~\ref{tab:trading_strategies1}--\ref{tab:trading_strategies2} and \ref{tab:trading_strategies3}--\ref{tab:trading_strategies11}, allowing us to evaluate how different activation functions influence neural network behavior in financial contexts. Our experiments are performed on an HPC cluster integrated with two NVIDIA H100 NVL GPUs, which significantly reduces training time and speeds up computation. This configuration ensures consistent GPU acceleration while providing efficient resource sharing across large computational workloads. 

\subsection{Performance Matrices for Trading Strategy}\label{4.2}
The performance of a trading strategy can vary significantly with different activation functions in a neural network, as each function influences the model’s ability to learn and generalize from data. To assess predictive accuracy and robustness, we evaluate our models using three key performance metrics: Root Mean Squared Error (RMSE), Mean Absolute Error (MAE), and the R-squared score (R\(^2\)). These metrics serve as essential indicators of the model’s reliability and precision.

Here, we investigate how different activation functions impact trading strategy performance across portfolios of varying sizes, specifically with \(10\) and \( 50\) assets, as shown in Tables~\ref{tab:trading_strategies1} and~\ref{tab:trading_strategies2}. For each portfolio size, we train multiple models using distinct activation functions and compare their performance across the chosen metrics. This comparative analysis reveals how the number of assets in the portfolio affects the sensitivity of each activation function in terms of predictive accuracy.

The variations in performance are summarized in Tables~\ref{tab:trading_strategies1} and~\ref{tab:trading_strategies2}, illustrating the effect of activation function choice on RMSE, MAE, and R\(^2\). This approach not only identifies the most effective activation function under different financial settings but also provides insights into how ridgelet-based neural network architecture interacts with portfolio complexity to influence overall model performance.

\begin{table}[ht!]
    \centering
     \renewcommand{\arraystretch}{5}
\rowcolors{2}{blue!10}{blue!10} 
    \renewcommand{\arraystretch}{1.2} 
    \begin{tabular}{|p{2cm}|p{2cm}|p{2cm}|p{2cm}|p{2cm}|p{2cm}|}
     \rowcolor{green!30} 
        \hline
        \textbf{Metric} & \textbf { SiLU }  &\textbf {GELU} & \textbf{ Mish} &\textbf {ReLU} &\textbf { HRDNN} \\
        \hline
        { RMSE} & {1.4922}   & {1.5416 } & 1.5384   & {1.7196} & 1.6637 \\
      { MAE} & {0.9094}  & {0.9378} & 0.9526 & {1.0926} & 1.0409 \\
       
        { R$^2$} & {0.9977} & {0.9977} & 0.9976 & {0.997} & 0.9972  \\
        \hline
    \end{tabular}
      \caption{Performance of trading strategies with various activation functions for 10 assets across key metrics.}

    \label{tab:trading_strategies1}
\end{table}

\begin{table}[h!]
    \centering
    \renewcommand{\arraystretch}{1.3}
\rowcolors{2}{blue!10}{blue!10} 
    \renewcommand{\arraystretch}{1.2} 
   \begin{tabular}{|p{2cm}|p{2cm}|p{2cm}|p{2cm}|p{2cm}|p{2cm}|}
     \rowcolor{green!30} 
        \hline
        \textbf{Metric} & \textbf { SiLU }  &\textbf {GELU} & \textbf{ Mish} &\textbf {ReLU} &\textbf { HRDNN} \\
        \hline
         { RMSE} & {1.8499}   & {1.9187 } &  2.2086  & {3.0853} & 2.5458 \\
      { MAE} & {1.1117}  & {1.1485} & 1.3882 & {2.2244} & 1.7278 \\
      
        { R$^2$} & {0.9983} & {0.9982} & 0.9976 & {0.9952} &  0.9967 \\
        \hline
    \end{tabular}
   \caption{Performance of trading strategies with various activation functions for 50 assets across key metrics.}
    \label{tab:trading_strategies2}
\end{table}

\subsection{Performance Ratios for Trading Strategy}\label{4.3}

Performance ratios play a crucial role in evaluating the profitability, risk profile, and overall efficiency of financial trading strategies. In this study, we focus on two widely used risk-adjusted performance metrics: the Sharpe ratio and the Sortino ratio. The Sharpe ratio measures the excess return relative to total volatility, offering a general assessment of risk-adjusted performance. In contrast, the sortino ratio refines this evaluation by considering only downside risk, making it especially relevant for strategies aimed at minimizing losses. Together, these metrics provide a comprehensive view of a strategy's ability to deliver returns under varying risk conditions.

Our analysis incorporates different transaction cost scenarios, including ZTC, PTC, and PSTC. These cost structures are used to assess the impact of activation function choices on trading strategy performance and to compare the outcomes against a standard benchmark, the buy-and-hold (B\&H) strategy. A summary of the results from this comparative study is presented in Tables~\ref{tab:trading_strategies3}--\ref{tab:trading_strategies11}.

\subsection{Comparison with the Linear Programming Approach}\label{4.4}

\begin{figure*}[!ht]
    \centering
    \includegraphics[width=14cm]{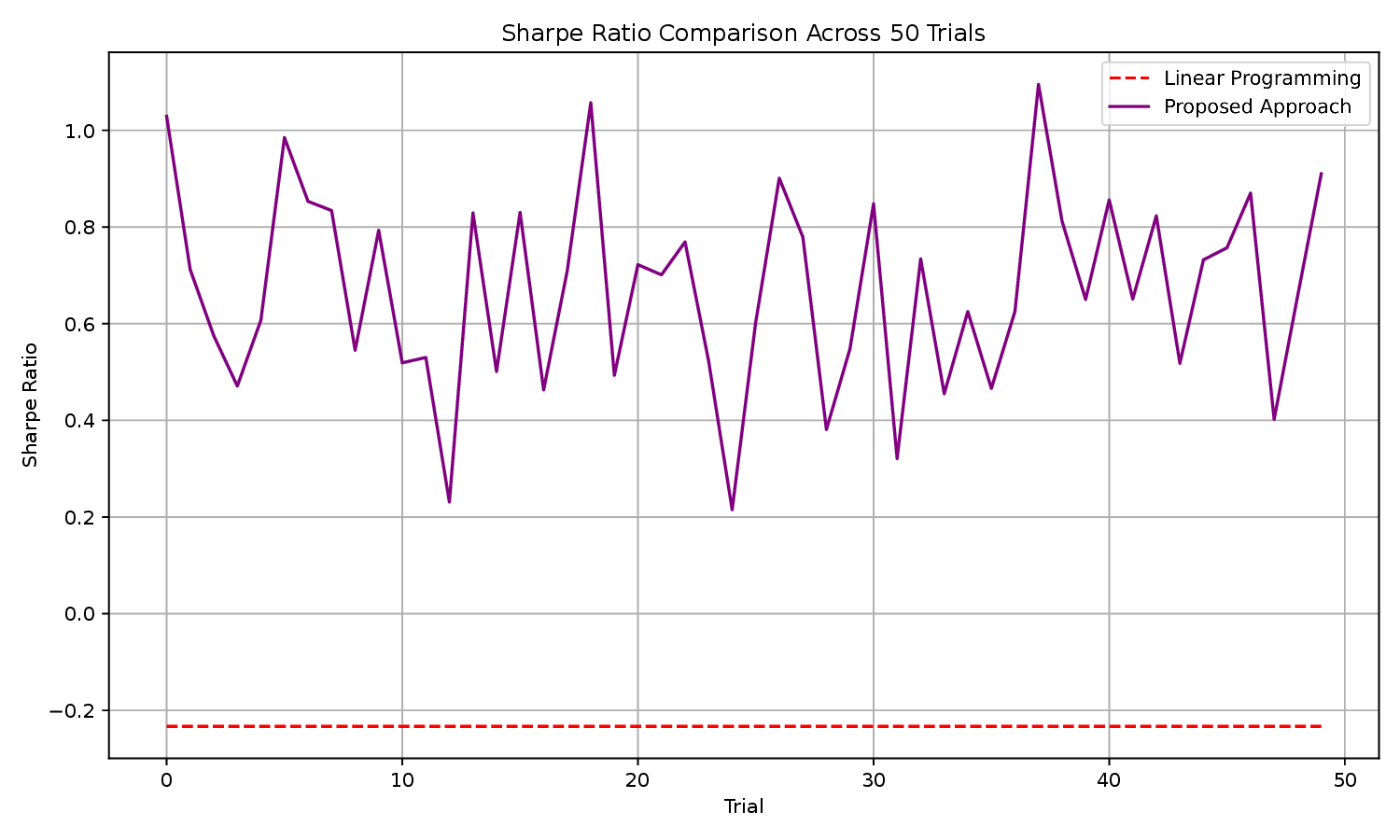}
    \caption{Comparison of Sharpe ratios between the Linear Programming approach and the proposed method across iterations (1--50).
}
    \label{fig:1}
\end{figure*}

In this section, we compare our proposed method with the linear programming (LP) approach introduced by Lütkebohmert. It is important to note that the LP method is computationally feasible only in low-dimensional scenarios, typically when the number of assets satisfies \( a \leq 3 \). To ensure consistency, we adopt the same experimental setup, focusing on the S\&P 500 index from the American stock market. The training period spans from January 1900 to December 2013, while the testing period covers September 2013 to December 2024. The trading strategy considers a one-month horizon with a single intermediate trade. Each strategy is trained across 1 to 50 iterations. We evaluate performance during the testing phase using the Sharpe ratio. Figure~\ref{fig:1} illustrates the Sharpe ratio as a function of the number of training iterations. Our results indicate that the proposed method consistently outperforms the LP-based strategy across most iterations. Notably, performance peaks around iteration 39, after which further training yields diminishing returns. This pattern suggests a trade-off between training time and generalization ability. Furthermore, the best performance measured by the highest Sharpe ratio is achieved when asset price constraints are relatively tight, indicating that adding stricter bounds does not hinder the identification of effective trading strategies. This finding highlights the robustness of our framework and its efficiency in navigating constrained financial environments.
\subsection{Trading in a Large Number of Stocks}\label{4.5}

This section evaluates the scalability and effectiveness of our trading strategy when applied to large, diverse portfolios. We begin with the case of 10 securities, selecting companies with the ticker symbols OKE, PG, RCL, SBUX, UNM, USB, VMC, WELL, WMB, and XOM. These stocks, all constituents of the S\&P 500 index, are chosen based on the Global Industry Classification Standard (GICS) to ensure sectoral diversification and to mitigate bias, particularly due to the dominance of IT firms in the U.S. market.

\begin{figure*}[!ht]
    \centering
    \includegraphics[width=15cm]{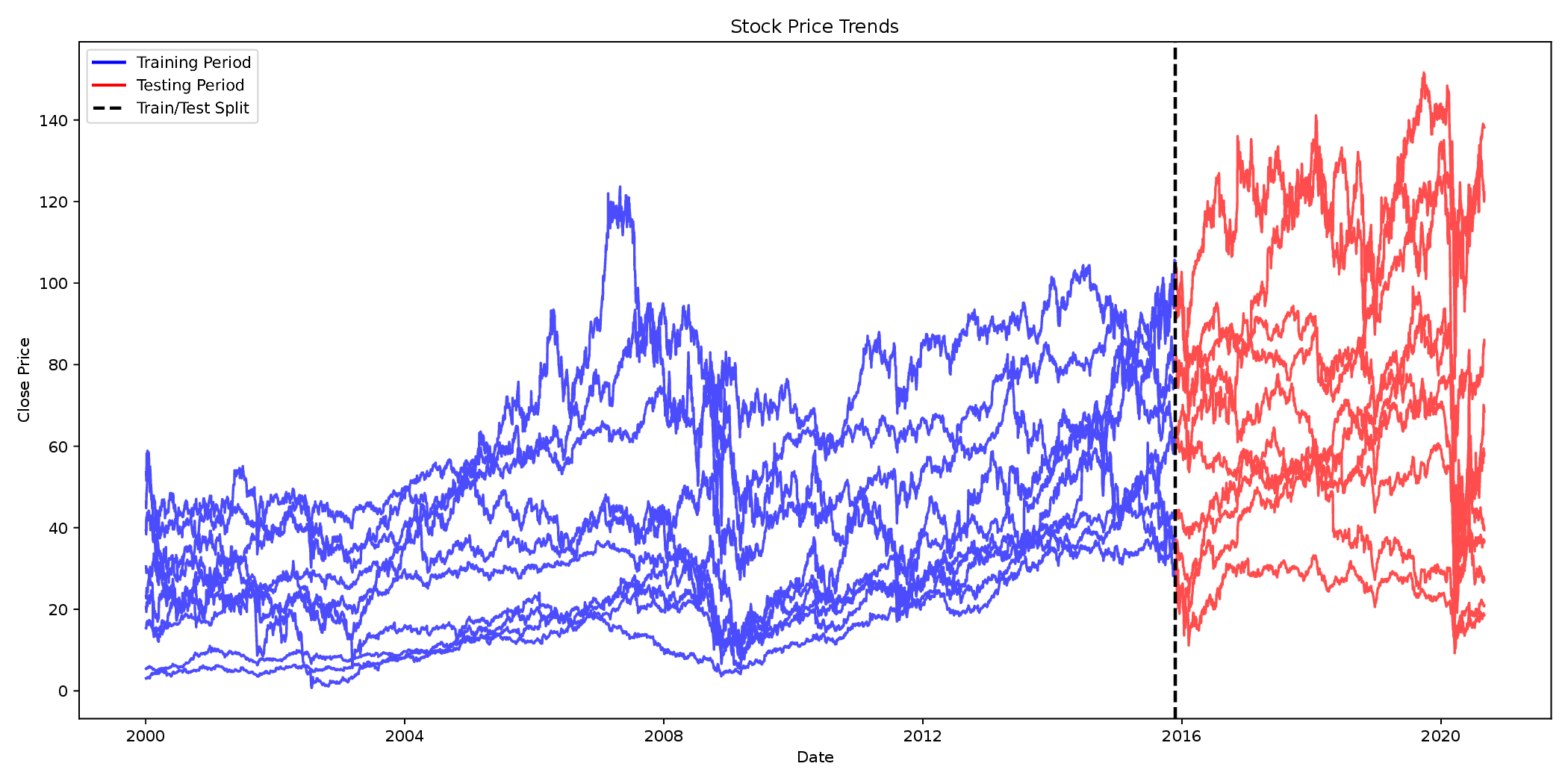}
    \caption{This figure illustrates the price evolution of 10 selected securities during the training period from 2000/01/03 to 2015/11/25 (shown in blue) and the testing period from 2015/11/26 to 2020/12/31 (shown in red).}
    \label{fig:2}
\end{figure*}

\begin{figure*}[!ht]
     \centering
    \includegraphics[width=15cm]{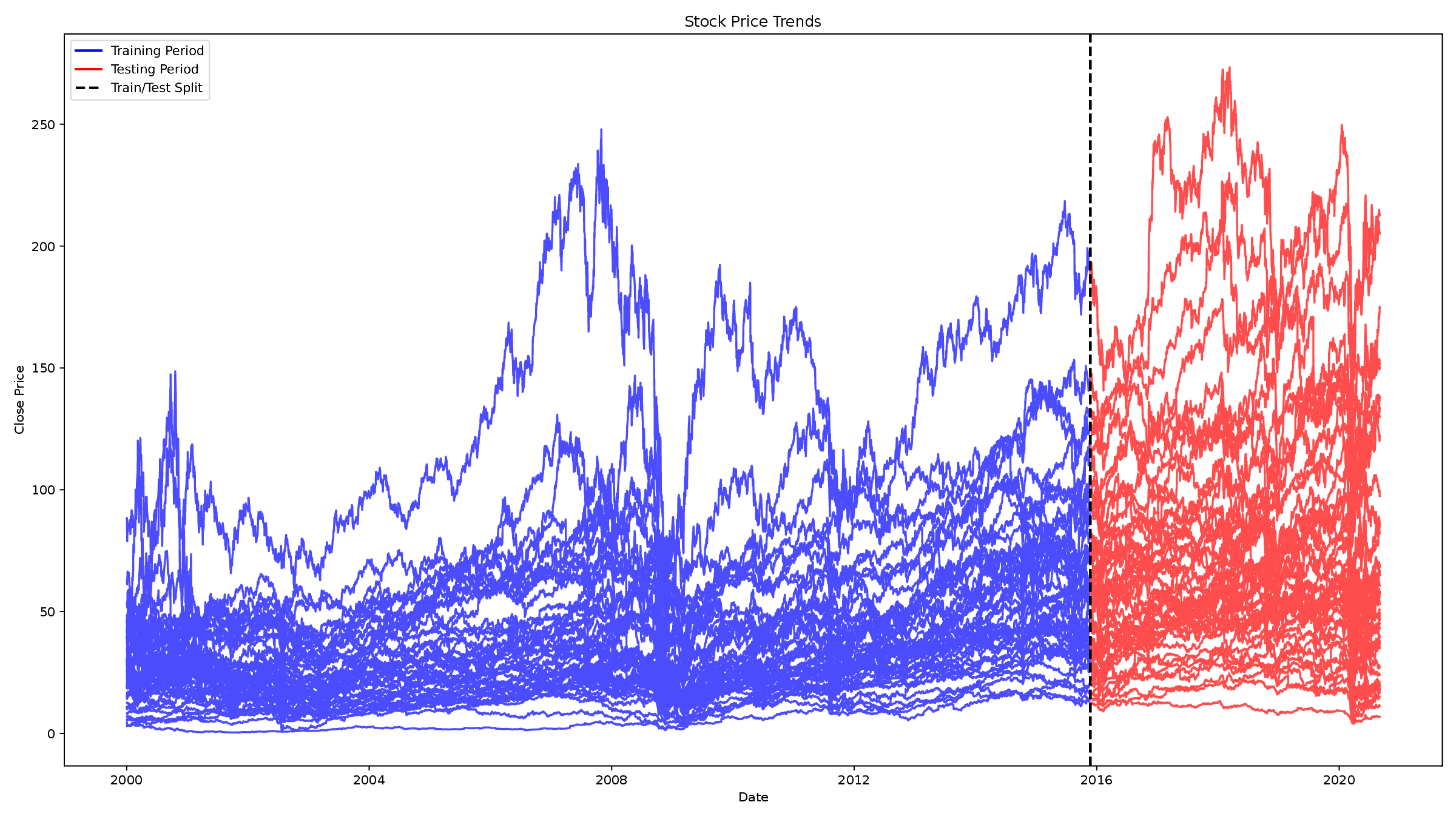}
  \caption{This figure illustrates the price evolution of 50 selected securities during the training period from 2000/01/03 to 2015/11/25 (shown in blue) and the testing period from 2015/11/26 to 2020/12/31 (shown in red).}
    \label{fig:6}
    \end{figure*}

The dataset is divided into a training period from January 3, 2000, to November 25, 2015, and a testing period from November 26, 2015, to December 31, 2020. Figures \ref{fig:2} and \ref{fig:6} illustrate the price evolution of the selected securities. Tables~\ref{tab:trading_strategies3}--\ref{tab:trading_strategies11} summarize the performance of the trained strategies during the testing phase. Compared to traditional methods such as linear programming and the B\&H strategy, our approach demonstrates superior handling of high-dimensional financial data and consistently outperforms the market across all scenarios considered.

Tables~\ref{tab:trading_strategies3}, \ref{tab:trading_strategies4}, and \ref{tab:trading_strategies5} provide a detailed comparison of trading strategies trained using different neural network activation functions, including SiLU, GELU, Mish, ReLU, and  HRDNN configuration, all applied to a portfolio of 10 securities. The evaluation is conducted under various transaction cost scenarios, assuming a frictionless market. As a benchmark, the B\&H strategy is also included. The models are assessed using multiple performance indicators such as total and average profits, percentage of profitable trades, maximum and minimum trade profits, and risk-adjusted measures including the Sharpe and Sortino ratios.

\begin{table*}[h!]
    \centering
     \renewcommand{\arraystretch}{1.3}
\rowcolors{2}{blue!10}{blue!10} 
    \renewcommand{\arraystretch}{1.2} 
    \resizebox{\textwidth}{!}{
     \begin{tabular}{l c c c c c c}
     \rowcolor{green!30} 
        \toprule
        \textbf{Zero Transaction Costs} & \textbf { SiLU } &\textbf {  B\&H} & \textbf {GELU} & \textbf{ Mish} &\textbf {ReLU} &\textbf { HRDNN} \\
        \midrule
        \textbf{Overall Profit} & {121.93} & -283.70 & {93.99 } & 111.22 & {113.34} & \textbf { 138.13} \\
        \textbf{Average Profit} & {1.02} & -2.36 & {0.78 } & 0.93 & {0.94} &\textbf {  1.15} \\
        \textbf{\% of Profitable Trades} & {54.22} & 60.0 & {55.22 } & 54.6  & {56.65} & \textbf { 58.48} \\
        \textbf{Maximum Profit} & {198.22} & 1006.58 & {193.72} & 221.23 & {221.97} &  \textbf { 226.58} \\
        \textbf{Minimum Profit} & {-149.24} & -2027.22 & {-188.94} & -218.5 & {-287.89} & \textbf { -283.1} \\
        \textbf{Sharpe Ratio} & {0.1643} & -0.043 & {0.082} & 0.1159  & {0.0979} &  \textbf { 0.1325} \\
            \textbf{Sortino Ratio} & {0.2659} & -0.037 & {0.1392 } & 0.1599 & {0.1333} &  \textbf { 0.1921} \\
        \bottomrule
   \end{tabular}
}
  \caption{Performance of trading strategies with different activation functions for 10 securities, compared under ZTC and B\&H.}
\label{tab:trading_strategies3}
\end{table*}
\begin{table*}[!ht]
    \centering
    \renewcommand{\arraystretch}{1.3}
\rowcolors{2}{blue!10}{blue!10} 
    \renewcommand{\arraystretch}{1.2} 
  \begin{tabular}{l c c c c c c}
     \rowcolor{green!30} 
        \toprule
        \textbf{Proportional Transaction Costs} & \textbf { SiLU } &\textbf {  B\&H} & \textbf {GELU} & \textbf{ Mish} &\textbf {ReLU} &\textbf { HRDNN} \\
        \midrule
        \textbf{Overall Profit} & \textbf {94.58} & -449.06 & {67.07} & 57.56 & {66.38 } &  { 75.23} \\
        \textbf{Average Profit} & {0.79} & -3.74 & {0.56} & 0.48  & {0.55} & \textbf { 0.63} \\
        \textbf{\% of Profitable Trades} & {52.12} & 60.0 & {53.37} & 50.72 & {55.15} & \textbf { 65.78} \\
        \textbf{Maximum Profit} & {197.13} & 1005.38 & {193.55} & 227.97  & {220.56} & \textbf { 265.02}\\
        \textbf{Minimum Profit} & {-149.24} & -2028.37 & {-190.33 } & -205.16 & {-288.94} &\textbf {  -385.89} \\
        \textbf{Sharpe Ratio} & {0.1085} & -0.068 & {0.0325} & 0.0405 & {0.0346} & \textbf { 0.0257} \\
        \textbf{Sortino Ratio} & {0.2029} & -0.058 & {0.0876} &  0.1526  & {0.0646} & \textbf { 0.0569} \\
        \bottomrule
    \end{tabular}
  \caption{Performance of trading strategies with different activation functions for 10 securities, compared under PTC  and B\&H.}
    \label{tab:trading_strategies4}
\end{table*}
\begin{table*}[!ht]
    \centering
     \renewcommand{\arraystretch}{1.3}
\rowcolors{2}{blue!10}{blue!10} 
    \renewcommand{\arraystretch}{1.2} 
  \begin{tabular}{l c c c c c c}
     \rowcolor{green!30} 
        \toprule
        \textbf{Per-Share Transaction Costs} & \textbf { SiLU } &\textbf {  B\&H} & \textbf {GELU} & \textbf{ Mish} &\textbf {ReLU} &\textbf { HRDNN} \\
        \midrule
        \textbf{Overall Profit} & {87.57} & -523.70 & {64.65} & 87.57 & {88.05 } & \textbf { 94.66} \\
        \textbf{Average Profit} & {0.73} & -4.36 & {0.54 } &0.73   & {0.73} &  \textbf { 0.79}\\
        \textbf{\% of Profitable Trades} & {52.18} & 60.0 & {52.9} & 52.18  & {56.07} &  \textbf { 56.92}\\
        \textbf{Maximum Profit} & {196.86} & 1004.58  & {201.4} & 196.86  & {237.01} &\textbf { 231.34} \\
        \textbf{Minimum Profit} & {-152.14} & -2029.22 & {-201.33} &-152.14   & {-316.92} & \textbf { -295.04} \\
        \textbf{Sharpe Ratio} & {0.0835} & -0.079 & {0.0265} & 0.0835 & {0.054} & \textbf { 0.0722} \\
        \textbf{Sortino Ratio} & {0.1623} & -0.068 & {0.0805} & 0.1623 & {0.080} & \textbf { 0.1059 }\\
        \bottomrule
     \end{tabular}
    \caption{Performance of trading strategies with different activation functions for 10 securities, compared under PSTC and B\&H.}
    \label{tab:trading_strategies5}
\end{table*}
\begin{table*}[h!]
    \centering
     \renewcommand{\arraystretch}{2}
\rowcolors{2}{blue!10}{blue!10} 
    \renewcommand{\arraystretch}{1.2} 
    \resizebox{\textwidth}{!}{
  \begin{tabular}{l c c c c c c}
       \rowcolor{green!30} 
        \toprule
        \textbf{Zero Transaction Costs} & \textbf { SiLU } &\textbf {  B\&H} & \textbf {GELU} & \textbf{ Mish} &\textbf {ReLU} &\textbf { HRDNN} \\
        \midrule
        \textbf{Overall Profit} & {103.18} & -331.27 & { 151.71} & 60.51  & {365.41} & \textbf { 489.4} \\        \textbf{Average Profit} & { 0.86} & -2.76 & {1.26  } &  0.5& {3.05 } &\textbf {  4.08}\\
        \textbf{\% of Profitable Trades} & { 61.9} & 61.67 & {63.1} & 62.22 & {62.87 } & \textbf { 64.88}\\
        \textbf{Maximum Profit} & {310.92} & 2289.62 & {371.87} & 388.56  & { 584.73} & \textbf { 688.74} \\
        \textbf{Minimum Profit} & {-433.46} & -3942.38 & {-564.1} & -641.5 & {-1114.42} & \textbf { -1459.01}  \\
        \textbf{Sharpe Ratio} & {0.0724} & -0.022 & {0.0887} & 0.0296 & {0.0857 } & \textbf { 0.1258 }\\
        \textbf{Sortino Ratio} & {0.0741} & -0.021 & {0.0822} & 0.0314 & {  0.0817 } & \textbf { 0.0974 }\\
        \bottomrule
     \end{tabular}
    }
   \caption{Performance of trading strategies with different activation functions for 20 securities, compared under ZTC and B\&H.}
    \label{tab:trading_strategies6}
\end{table*}

\begin{table*}[!ht]
    \centering
         \renewcommand{\arraystretch}{1.3}
\rowcolors{2}{blue!10}{blue!10} 
    \renewcommand{\arraystretch}{1.2} 
 \begin{tabular}{l c c c c c c}
    \rowcolor{green!30} 
        \toprule
        \textbf{Proportional Transaction Costs} & \textbf{SiLU} & \textbf {  B\&H} & \textbf{GELU} & \textbf { Mish} & \textbf{ReLU} & \textbf {HRDNN} \\
        \midrule
        \textbf{Overall Profit} & {371.79  } & -839.85 & {452.54 } & 380.84 & {599.74} &\textbf { 871.49}  \\
        \textbf{Average Profit} & {3.1} &  -7.0& {3.77 } & 3.17 & {5.0} & \textbf { 7.26}  \\
        \textbf{\% of Profitable Trades} & { 60.9} & 61.67  & {61.82 } & 60.77  & { 60.02} & \textbf { 64.13}  \\
        \textbf{Maximum Profit} & {514.85} & 2999.84 & {624.71} &685.06  & {921.45 } &  \textbf { 1150.14 }\\
        \textbf{Minimum Profit} & {-628.06 } & -5055.84 & {-869.24 } &-925.4  & {-1703.46 } &\textbf { -2399.37}  \\
        \textbf{Sharpe Ratio} & {0.1797} & -0.044& {0.1725} &0.1225 & { 0.052 } & \textbf { 0.1313} \\
        \textbf{Sortino Ratio} & {0.174} & -0.041 & {0.1556} & 0.1156 & {0.0379} & \textbf { 0.1045} \\
        \bottomrule
     \end{tabular}
   \caption{Performance of trading strategies with different activation functions for 30 securities, compared under PTC and B\&H.}

    \label{tab:trading_strategies7}
\end{table*}

\begin{table*}[!ht]
    \centering
      \renewcommand{\arraystretch}{1.3}
\rowcolors{2}{blue!10}{blue!10} 
    \renewcommand{\arraystretch}{1.2} 
  \begin{tabular}{l c c c c c c}
    \rowcolor{green!30} 
        \toprule
        \textbf{Per-Share Transaction Costs} & \textbf{SiLU} & \textbf {  B\&H} & \textbf{GELU} &\textbf { Mish} & \textbf{ReLU} & \textbf { HRDNN} \\
        \midrule
        \textbf{Overall Profit} & {637.35} & -1475.84 & {  757.57 } & 613.16  & {977.64 } & \textbf {  1042.54}\\
        \textbf{Average Profit} & {5.31} & -12.3 & {6.31} & 5.11 & {8.15} & \textbf { 8.69 }  \\
        \textbf{\% of Profitable Trades} & {62.32 } & 61.67 & {64.28 } & 61.43  & {64.97} &\textbf {  66.82} \\
        \textbf{Maximum Profit} & { 861.13} & 3599.23 & { 999.59} & 1043.46 & {1490.76} & \textbf { 1382.5 } \\
        \textbf{Minimum Profit} & {-880.71} & -6093.67 & {-1284.36} &-1316.46   & {-3305.51} & \textbf { -2671.33} \\
        \textbf{Sharpe Ratio} & {0.1928} & -0.063 & {0.1824} & 0.131 & {0.1064} & \textbf { 0.1372 }  \\
        \textbf{Sortino Ratio} & { 0.1899} & -0.058  & {0.1607} & 0.1214 & {0.0826 } & \textbf { 0.1057}\\
        \bottomrule
   \end{tabular}
   \caption{Performance of trading strategies with different activation functions for 40 securities, compared under PSTC and B\&H.}

    \label{tab:trading_strategies8}
\end{table*}

\begin{table*}[!ht]
    \centering
          \renewcommand{\arraystretch}{1.3}
\rowcolors{2}{blue!10}{blue!10} 
    \renewcommand{\arraystretch}{1.2} 
\begin{tabular}{l c c c c c c}
     \rowcolor{green!30} 
        \toprule
        \textbf{Zero Transaction Costs} & \textbf{SiLU} & \textbf {  B\&H} & \textbf{GELU} &\textbf { Mish} & \textbf{ReLU} & \textbf { HRDNN} \\
        \midrule
        \textbf{Overall Profit} & {948.63} & 327.79 & { 727.7 } & -523.26 & {636.73} & \textbf {1548.5} \\
        \textbf{Average Profit} & { 7.91} & 2.73 & {6.06} & -4.36 & {5.31} & \textbf { 12.9 }\\
        \textbf{\% of Profitable Trades} & {64.6} & 65.0 & {61.9} & 65.0 & {62.78 } & \textbf { 67.0} \\
        \textbf{Maximum Profit} & {1103.61} & 3991.14 & {1101.62 } & 3984.52 & {1198.78} & \textbf { 1639.17} \\
        \textbf{Minimum Profit} & {-958.17} & -7118.32 & {-978.98} & -7124.48 & {-1174.59} & \textbf { -3120.32 }\\
        \textbf{Sharpe Ratio} & { 0.2562} & 0.012 & {0.1947} & -0.019 & {0.1452} & \textbf { 0.1777} \\
        \textbf{Sortino Ratio} & {0.2623 } & 0.011 & {0.2049} & -0.017 & { 0.1416} & \textbf {0.1377} \\
        \bottomrule
   \end{tabular}
   \caption{Performance of trading strategies with different activation functions for 50 securities, compared under ZTC and B\&H.}
    \label{tab:trading_strategies9}
\end{table*}

\begin{table*}[!ht]
    \centering
            \renewcommand{\arraystretch}{1.3}
\rowcolors{2}{blue!10}{blue!10} 
    \renewcommand{\arraystretch}{1.2} 
    \begin{tabular}{l c c c c c c}
     \rowcolor{green!30} 
        \toprule
        \textbf{Proportional Transaction Costs} & \textbf{SiLU} &\textbf{ B\&H} & \textbf{GELU} & \textbf{Mish} & \textbf{ReLU} &\textbf{HRDNN} \\
        \midrule
        \textbf{Overall Profit} & { 727.7} & -523.26 & { 819.92 } &577.49  & {868.7 } & \textbf{1068.31 }\\
        \textbf{Average Profit} & { 6.06} & -4.36 & {6.83} & 4.81 & {7.24 } & \textbf{8.9} \\
        \textbf{\% of Profitable Trades} & {61.9} & 65.0 & {63.53} & 61.58  & {63.55 } & \textbf{65.38} \\
        \textbf{Maximum Profit} & {1101.62} & 3984.52 & { 1218.64} &1278.04  & {1777.29 } &\textbf{1610.55}  \\
        \textbf{Minimum Profit} & {-978.98} & -7124.48 & {-1419.34} & -1465.63  & {-4018.29 } & \textbf{-3126.47} \\
        \textbf{Sharpe Ratio} & { 0.1947} &  -0.019 & {0.1801} & 0.1142 & {0.0781} &\textbf{0.1206}  \\
        \textbf{Sortino Ratio} & {0.2049 } & -0.017 & {0.1667} & 0.1166 & {0.0627 } & \textbf{0.0942 }\\
        \bottomrule
   \end{tabular}
      \caption{Performance of trading strategies with different activation functions for 50 securities, compared under PTC and B\&H.}

    \label{tab:trading_strategies10}
\end{table*}

\begin{table*}[!ht]
    \centering
    \rowcolors{2}{blue!10}{blue!10} 
    \renewcommand{\arraystretch}{1.2} 
    
  \begin{tabular}{l c c c c c c}
        \rowcolor{green!30} 
        \toprule
        \textbf{Per-Share Transaction Costs} & \textbf{SiLU} & \textbf{ B\&H} & \textbf{GELU} & \textbf{Mish} & \textbf{ReLU} &\textbf{HRDNN} \\
        \midrule
        \textbf{Overall Profit} & {636.73} & -872.20 & { 783.52 } &  450.06 & {934.47} &\textbf{1036.16 } \\
        \textbf{Average Profit} & { 5.31} & -7.27 & {6.53} & 3.75 & {7.79} &  \textbf{8.63}\\
        \textbf{\% of Profitable Trades} & {62.78 } &65.0 & {64.03} &62.5   & {64.47} & \textbf{66.23 } \\
        \textbf{Maximum Profit} & {1198.78} & 3981.144 & {1366.01 } & 1387.04 & {1925.53} &\textbf{1829.21}  \\
        \textbf{Minimum Profit} & {-1174.59} & -7128.32 & {-1742.62} & -1772.34 & {-4577.7} & \textbf{-3745.27} \\
        \textbf{Sharpe Ratio} & { 0.1452} & -0.032 & {0.1419} &0.0735  & {0.0738 } & \textbf{0.099}  \\
        \textbf{Sortino Ratio} & { 0.1416} & -0.028 & {0.1232} & 0.0669 & {0.0573 } &\textbf{0.0756}   \\
        \bottomrule
    \end{tabular}
   \caption{Performance of trading strategies with different activation functions for 50 securities, compared under PSTC and B\&H.}
\label{tab:trading_strategies11}
\end{table*}


To further assess scalability and robustness, we extend our experiments to larger portfolios with 20, 30, 40, and 50  securities, as shown in Tables~\ref{tab:trading_strategies6}, \ref{tab:trading_strategies7}, \ref{tab:trading_strategies8}, \ref{tab:trading_strategies9}, \ref{tab:trading_strategies10}, and  \ref{tab:trading_strategies11}.  Unlike traditional techniques that are restricted to low-dimensional settings, our model demonstrates flexibility across varying portfolio sizes. We observe that the choice of activation function significantly affects performance metrics across these dimensions. Notably, HRDNN consistently outperforms individual functions, highlighting its ability to model complex market patterns and enhance trading outcomes.

Summary of Tables~\ref{tab:trading_strategies3}--\ref{tab:trading_strategies11}, the B\&H strategy performs consistently poorly across all transaction cost levels (ZTC, PTC, and PSTC), generating significant losses and negative Sharpe and Sortino ratios.
In terms of activation functions, HRDNN performs the best, generating the largest average and total profits under ZTC, while maintaining competitiveness under PTC, and PSTC. SiLU ranks second overall and offers consistent profitability across scenarios, whereas Mish and ReLU exhibit moderate gains with occasionally greater maximum profits. With lower returns and risk-adjusted measures, GELU trails behind. B\&H continues to be the weakest baseline, whereas HRDNN is the most reliable option overall, combining stability and profitability. Our model is consistently superior to the B\&H approach when we extend it to higher dimensions. Our model generates robust gains and favorable risk-adjusted ratios across 20, 30, 40, and 50 securities across various transaction cost settings, but B\&H continues to produce losses. HRDNN is the most reliable function, followed by ReLU and SiLU.

Figures~\ref{fig:2}–\ref{fig:6} illustrate the price evolution of portfolios containing 10 to 50 securities, covering a wide range of industries and volatility patterns. As the number of assets increases, the portfolios become more diversified, offering richer information and broader trading opportunities. Our model adapts to these heterogeneous dynamics by capturing cross-asset dependencies, whereas the B\&H strategy merely averages across assets and remains exposed to extended drawdowns. In higher-dimensional settings, this diversification effect amplifies the model’s advantage, leading to more consistent and profitable outcomes while B\&H continues to underperform.

\section{Conclusion }\label{5}
In this paper, we propose a novel framework that integrates deep neural networks with Ridgelet Transform techniques to design profitable trading strategies under model uncertainty. Ridgelet Transform enhances the expressive power of neural networks by enabling effective approximation of complex, high-dimensional payoff structures that are central to arbitrage strategies. In particular, we prove the convergence of a penalized functional that incorporates HRDNN approaches towards the minimal conditional super-replication price of profitable trading strategies. Theoretically, this approach allows HRDNN to be optimized in order to determine profitable trading strategies.  We also emphasized the importance of activation functions in achieving robust approximation and adaptability in uncertain environments. Our systematic evaluation showed that the choice of activation function significantly influences both predictive accuracy and profitability. Our research shows that the HRDNN generates higher returns and performs significantly better than the traditional B\&H approach. To validate the practicality of our approach, we conducted empirical studies, which demonstrated consistently profitable outcomes, even during periods of high market volatility. This work can be further extended to incorporate diverse asset classes, generative modeling approaches, adaptive activation mechanisms, and validation under real-time market conditions.

\bibliographystyle{siamplain}
\bibliography{references}
\end{document}